\documentclass[a4paper,10pt]{article}
\usepackage{latexsym}
\usepackage{amsmath}
\usepackage{amssymb}
\usepackage{amscd}
\usepackage{graphicx}

\newtheorem{Th}{Theorem}[section]
\newtheorem{Cor}{Corollary}[section]
\newtheorem{Prop}{Proposition}[section]
\newtheorem{Lem}{Lemma}[section]
\newtheorem{Def}{Definition}[section]
\newtheorem{Rem}{Remark}[section]
\newtheorem{Ex}{Example}[section]

\newcommand{\bet}{\begin{Th}}
\newcommand{\ent}{\stepcounter{Cor}
   \stepcounter{Prop}\stepcounter{Lem}\stepcounter{Def}
   \stepcounter{Rem}\stepcounter{Ex}\end{Th}}

\newcommand{\bec}{\begin{Cor}}
\newcommand{\enc}{\stepcounter{Th}
   \stepcounter{Prop}\stepcounter{Lem}\stepcounter{Def}
   \stepcounter{Rem}\stepcounter{Ex}\end{Cor}}
\newcommand{\bep}{\begin{Prop}}
\newcommand{\enp}{\stepcounter{Th}
   \stepcounter{Cor}\stepcounter{Lem}\stepcounter{Def}
   \stepcounter{Rem}\stepcounter{Ex}\end{Prop}}
\newcommand{\bel}{\begin{Lem}}
\newcommand{\enl}{\stepcounter{Th}
   \stepcounter{Cor}\stepcounter{Prop}\stepcounter{Def}
   \stepcounter{Rem}\stepcounter{Ex}\end{Lem}}
\newcommand{\bef}{\begin{Def}}
\newcommand{\enf}{\stepcounter{Th}
   \stepcounter{Cor}\stepcounter{Prop}\stepcounter{Lem}
   \stepcounter{Rem}\stepcounter{Ex}\end{Def}}
\newcommand{\ber}{\begin{Rem}}
\newcommand{\enr}{
   \stepcounter{Th}\stepcounter{Cor}\stepcounter{Prop}
   \stepcounter{Lem}\stepcounter{Def}\stepcounter{Ex}\end{Rem}}
\newcommand{\bee}{\begin{Ex}}
\newcommand{\ene}{
   \stepcounter{Th}\stepcounter{Cor}\stepcounter{Prop}
   \stepcounter{Lem}\stepcounter{Def}\stepcounter{Rem}\end{Ex}}
\newcommand{\Proof}{\noindent{\it Proof\,}:\ }

\newcommand{\R}{\mathbf{R}}
\newcommand{\C}{\mathbf{C}}

\newcommand{\lon}{\longrightarrow}

\newcommand{\pa}{\partial}

\newcommand{\QED}{\hfill$\Box$\par}

\begin{document}

\title{Singularities of improper affine spheres and \\ 
surfaces of constant Gaussian curvature} 

\author{Go-o ISHIKAWA\thanks{Partially supported by 
Grants-in-Aid for Scientific Research, No. 14340020. } \ 
and Yoshinori MACHIDA\thanks{Partially supported by 
Grants-in-Aid for Scientific Research, (C) No. 14540097. }}
\date{ }
\renewcommand{\thefootnote}{\fnsymbol{footnote}}
\footnotetext{
{\hspace{-0.54truecm}}2000 {\it Mathematics Subject Classification}:  
Primary. 
53C42, 
53A15, 
Secondary. 
58K40,  
58A15. 
}
\maketitle

{\small 
\noindent {\it Abstract}: 
We study the equation for improper (parabolic) affine spheres 
from the view point of contact geometry 
and provide the generic classification of singularities 
appearing in geometric solutions to the equation as well as their duals. 
We also show the results for 
surfaces of constant Gaussian curvature
and for developable surfaces. 
In particular we confirm that generic singularities appearing in such a surface 
are just cuspidal edges and swallowtails.  
}

\

\begin{center}
{\it Dedicated to Professor Tohru Morimoto on his 60th birthday}
\end{center}

\section{Introduction.}

Let $f(x, y)$ be a $C^\infty$ 
function on $\R^2$ 
satisfying the unimodular Hessian equation: 
$$
{\mathrm{Hess}}(f) = 
\left\vert 
\begin{array}{cc}
\dfrac{\pa^2 f}{\pa x^2} & \dfrac{\pa^2 f}{\pa x\pa y}\\
\dfrac{\pa^2 f}{\pa y\pa x} & \dfrac{\pa^2 f}{\pa y^2}
\end{array}
\right\vert
= \pm 1. 
$$
Then the graph 
$z = f(x, y)$ in $\R^3$ is an {\it improper 
affine sphere}  
with the affine normal vector field $\pa/\pa z$ \cite{NS}. 
In this note we study the equations \lq\lq Hessian $= \pm 1$" and 
singularities of improper affine spheres. 
Also we study 
the equations of constant Gaussian curvature 
$K = c$ 
for surfaces in $\R^3$. 
We 
provide the results on singularities 
of geometric solutions to $K = -1$ (\lq\lq pseudo-spherical surfaces"), 
$K = 1$ (\lq\lq sphere-like surfaces") and $K = 0$ (\lq\lq developable surfaces"). 

The importance of the study on singularities 
comes from the well-known classical results: 
A smooth global solution in $\R^3$ to ${\mathrm{Hess}} = 1$ 
is the graph of a quadratic polynomial function (J\"{o}rgens' theorem
\cite{J}\cite{C}\cite{P}, cf.{\,}Bernstein's theorem \cite{B2}\cite{N}). 
Therefore other solutions necessarily have singularities. 
Besides we know: A compact surface in $\R^3$ with constant positive 
Gaussian curvature is a sphere (Liebmann's theorem\cite{O}). 
Therefore compact solution surfaces to 
$K = 1$ other than spheres have necessarily singularities. 
Moreover, 
there are no complete surface in $\R^3$ with constant negative 
Gaussian curvature (Hilbert's theorem\cite{Hil}, cf.{\,}\cite{Ts}). 
Therefore solution surfaces to the equation $K = -1$ have necessarily singularities. 

As well-known, the equations 
${\mathrm{Hess}} = c$ and 
$K = c$ are regarded as Monge-Amp\`{e}re 
equations, and they have been studied from both geometric and analytic 
aspects. Note that, for a surface $z = f(x, y)$, we have 
$K = \dfrac{f_{xx}f_{yy}-f_{xy}^2}{(1+f_x^2+f_y^2)^2}$ and 
the equation $K = c$ turns to be 
$$
f_{xx}f_{yy}-f_{xy}^2 = c(1+f_x^2+f_y^2)^2. 
$$
Then we observe that 
there are associated to these equations an additional 
geometric structure, the structure of {\it Lagrangian pair} (\cite{IMa}), and 
based on that, we proceed to detailed study on the singularities of solutions 
beyond usually expected. 

\ 

Consider the equation $f_{xx}f_{yy} - f_{xy}^2 = c$ 
for a surface $z = f(x, y)$ in $xyz$-space $\R^3$.  
Geometrically the equation can be written into the differential system
$$
\omega = c dx \wedge dy - dp \wedge dq = 0, \quad 
\theta = dz - pdx - qdy = 0, 
$$
on $xyzpq$-space $\R^5$, $p, q$ representing $z_x, z_y$ respectively. 
Then we have the contact distribution $D = \{ \theta = 0 \}$ 
with the symplectic structure $d\theta\vert_D$ 
in the tangent bundle $T\R^5$. 
Moreover 
$E_1 = \{ v \in D \mid i_v(dx\wedge dy) = 0 \}$ and 
$E_2 = \{ v \in D \mid i_v(dp\wedge dq) = 0 \}$ 
are Lagrangian subbundles of $D$, where $i_v$ 
denotes the interior product by $v$. In fact 
$E_1$ (resp. $E_2$) is generated by 
$\dfrac{\pa}{\pa p}, \dfrac{\pa}{\pa q}$ 
(resp. by $\dfrac{\pa}{\pa x} + p\dfrac{\pa}{\pa z}, 
\dfrac{\pa}{\pa y} + q\dfrac{\pa}{\pa z}$ ). 
Then the double Legendrian fibration 
$$\begin{CD}
     & @. & M = \R^5 & @. & \\
     \pi_1 & \swarrow & @. & \searrow & \pi_2 \\
     W_1 = \R^3 & @. & &  & W_2 = \R^3 
\end{CD}$$
is induced. 
The first projection $\pi_1$ is defined by $(x, y, z, p, q) \mapsto (x, y, z)$ 
and the second projection $\pi_2$ is defined by 
$(x, y, z, p, q) \mapsto (p, q, px + qy - z)$. 

\

The differential system associated to the equation $K = c$ in the 
Euclidean $3$-space 
$\R^3$ is defined on the unit tangent bundle $\R^3\times S^2$
of $\R^3$ by the $2$-form 
$$
\omega  = 
c(y_1dx_2\wedge dx_3 + y_2dx_3\wedge dx_1 + y_3dx_1\wedge dx_2) 
 - (y_1dy_2 \wedge dy_3 + y_2dy_3 \wedge dy_1 + y_3dy_1 \wedge dy_2),
$$
and the contact form 
$
\theta = y_1dx_1 + y_2dx_2 + y_3dx_3. 
$
Here $(x_1, x_2, x_3; y_1, y_2, y_3)$ is the system of coordinates on 
$\R^3\times \R^3$ restricted to $\R^3\times S^2$. 
We set $D = \{ \theta = 0\} \subset T(\R^3\times S^2)$, and 
two Lagrangian subbundles of $D$: 
$$
\begin{array}{rcl}
E_1 & = & \{ v \in D \mid i_v(y_1dx_2\wedge dx_3 + y_2dx_3\wedge dx_1 + y_3dx_1\wedge dx_2) = 0\} \\ 
 & = & \left.\left\{ v = \eta_1\dfrac{\pa}{\pa y_1} + \eta_2\dfrac{\pa}{\pa y_2} + \eta_3\dfrac{\pa}{\pa y_3}  \ \right\vert \ 
v  {\mbox{\rm{\ is\ tangent\ to\ }}} S^2 \right\}, 
\end{array}
$$
$$
\begin{array}{rcl}
E_2 & = & \{ v \in D \mid 
i_v(y_1dy_2 \wedge dy_3 + y_2dy_3 \wedge dy_1 + y_3dy_1 \wedge dy_2) = 0 \} \\ 
 & = & \left.\left\{ v = \xi_1\dfrac{\pa}{\pa x_1} + \xi_2\dfrac{\pa}{\pa x_2} + 
 \xi_3\dfrac{\pa}{\pa x_3} \ \right\vert \  \xi_1y_1 + \xi_2y_2 + \xi_3y_3 = 0 \right\}. 
\end{array}
$$
Then we have the double Legendrian fibration
$$\begin{CD}
     & @. & M = \R^3\times S^2 & @. & \\
     \pi_1 & \swarrow & @.  & \searrow & \pi_2 \\
     W_1 = \R^3 & @. & & & W_2 = \R\times S^2 , 
\end{CD}$$
\\
where $\pi_1(x, y) = x, \pi_2(x, y) = (x\cdot y, y)$ for $(x, y) \in \R^3\times S^2 
\subset \R^3\times \R^3$, 
and $x\cdot y$ is the inner product (the height function). 

\ 

In general, let $M$ be a contact manifold of dimension $2n+1$ with
a contact structure $D \subset TM$.
A {\it Monge-Amp\`{e}re system}
on $M$ 
is an exterior differential system ${\mathcal M}$ 
generated locally by a contact form $\theta$ for $D$ and an
$n$-form $\omega$ on $M$. 
The geometric formulation of Monge-Amp\`{e}re systems,
originally 
due to T. Morimoto \cite{M1}\cite{M2} and
V.V. Lychagin \cite{L}, 
naturally and intrinsically generalizes the classical Monge-Amp\`{e}re
equations, and describes several fundamental
equations in geometry, analysis and physics.
As explained above, we have Monge-Amp\`{e}re systems 
associated to the equations Hessian $= c$ and $K = c$ respectively. 

An immersion $f : N^{n} \to M$ from a manifold
$N$ of dimension $n$ is called a {\it Legendrian immersion} 
if $f_*(TN) \subset D (\subset TM)$, where $f_* : TN \to TM$ 
is the differential of $f$. 
A Legendrian immersion $f : N \to M$ 
is called a {\it geometric solution}
of a Monge-Amp\`{e}re system ${\mathcal M}$ if $f^*{\mathcal M} = 0$.
If ${\mathcal M}$ is given by a contact form $\theta$ and
an $n$-form $\omega$, 
then the condition reads $f^*\theta = 0, f^*\omega = 0$. 

\ 

We lift any solution surface $z = f(x, y)$ to ${\mathrm{Hess}} = c$ 
uniquely to a geometric solution to ${\mathrm{Hess}} = c$ in $\R^5$, 
by setting $p = f_x, q = f_y$. 
Then the $\pi_2$-projection of the lifting is nothing but the {\it affine dual} 
of the original surface. 
Also we lift any surface in $\R^3$ with $K = c$ with a fixed co-orientation, 
uniquely to a geometric solution to $K = c$ in $\R^3\times S^2$, 
by using the Gauss map (unit normals). 
The $\pi_2$-projection of the lifting is the {\it pedal surface} in 
$\R \times S^2 \cong \R^3 \setminus \{ 0 \}$ 
of the original surface. 
To describe singularities of solution surfaces, we start to study 
geometric solutions and then 
the both $\pi_1$ and 
$\pi_2$-projections of them. 

\ 

In this paper we show the following result: 

\bet
\label{cusp and swallowtail}
A generic geometric solution to ${\mathrm{Hess}} = 1$ 
{\rm (}resp. to ${\mathrm{Hess}} = -1$, $K = 1$, or $K = -1${\rm )} 
has only cuspidal edges and swallowtails as singularities. 
More  strictly, any geometric solution $N^2 \to M^5 = \R^5$ 
or $\R^3\times S^2$ 
to ${\mathrm{Hess}} = 1$ 
{\rm (}resp.{\,}to ${\mathrm{Hess}} = -1$, $K = 1$, or $K = -1${\rm )} 
can be locally approximated  near each point in $N$ in $C^\infty$ topology by 
a geometric solution $f : U \to M$ 
such that, for any $x_0 \in U$, 
one of the following assertions {\rm (i), (ii), (iii), (iv)} holds, 
with respect to 
the Legendrian fibrations $\pi_1 : M \to W_1 = \R^3$ and 
$\pi_2 : M \to W_2 = \R^3$ or $\R\times S^2$: 

{\rm (i)} 
$\pi_1\circ f : (U, x_0) \to W_1$ is an immersion at $x_0$, and 
$\pi_2\circ f : (U, x_0) \to W_2$ 
is an immersion at $x_0$.  

{\rm (ii)} 
$\pi_1\circ f$ has the cuspidal edge at $x_0$, and 
$\pi_2\circ f$ has the cuspidal edge at $x_0$.  

{\rm (iii)} 
$\pi_1\circ f$ has the swallowtail at $x_0$, and 
$\pi_2\circ f$ has the cuspidal edge at $x_0$.  

{\rm (iv)} 
$\pi_1\circ f$ has the cuspidal edge at $x_0$, and 
$\pi_2\circ f$ has the swallowtail at $x_0$.  

\ent

The significance of Theorem \ref{cusp and swallowtail} 
is twofold. 

First, 
Theorem \ref{cusp and swallowtail} is a collection of four theorems: 
Four results are independent to each other, 
since four equations have different properties to each other 
and we need to analyze geometric solutions for each equations. 
Neverthless we get the same list of generic singularities as a result. 

Note also that 
the equation ${\mathrm{Hess}} = c (c > 0)$ (resp. 
${\mathrm{Hess}} = c (c < 0)$) 
is equivalent to ${\mathrm{Hess}} = 1$ (resp. ${\mathrm{Hess}} = -1$), 
by the contactomorphism $(x, y, z, p, q) \mapsto 
(\sqrt{\vert c\vert}x, \sqrt{\vert c\vert}y, z, \dfrac{1}{\sqrt{\vert c\vert}}p, 
\dfrac{1}{\sqrt{\vert c\vert}}q)$. 
Similarly the equation $K = c (c > 0)$ (resp. $K = c (c < 0))$ 
is equivalent to $K = 1$ (resp. $K = -1)$ by the contactomorphism 
$(x, y) \to (\sqrt{\vert c\vert}x, y)$.

Second, the classification result of Theorem \ref{cusp and swallowtail} differs 
from that for generic Legendrian submanifolds: 

\bep
\label{general Legendre}
For a generic Legendrian immersion $N^2 \to \R^5$ {\rm (}not necessarily a 
geometric solution{\rm )}, 
one of the following holds: 

{\rm (a)}
$\pi_1\circ f$ is an immersion, and $\pi_2\circ f$ is an immersion at $x_0$. 

{\rm (b)}
$\pi_1\circ f$ is the cuspidal edge, and $\pi_2\circ f$ is an immersion at $x_0$. 

{\rm (c)}
$\pi_1\circ f$ is an immersion, and $\pi_2\circ f$ is the cuspidal edge 
at $x_0$. 

{\rm (d)}
$\pi_1\circ f$ is the cuspidal edge, and $\pi_2\circ f$ is the cuspidal edge 
at $x_0$.  

{\rm (e)}
$\pi_1\circ f$ is the swallowtail, and $\pi_2\circ f$ is an immersion at 
$x_0$. 

{\rm (f)}
$\pi_1\circ f$ is an immersion and $\pi_2\circ f$ is the swallowtail at 
$x_0$.  

\enp

Note that 
Proposition \ref{general Legendre} is a straightforward result from 
ordinary Legendrian singularity theory (\cite{AGV}\cite{Z}). 
In fact the result is described by 
the generic combination of two stratifications on the plane 
by $A_3, A_2, A_1$-singularities 
via two Legendrian fibrations $\pi_1$ and $\pi_2$ respectively. 

Contrary to Proposition \ref{general Legendre}, 
we observe the simultaneous occurrence of singularities 
of $\pi_1$ and $\pi_2$-projections generically, in Theorem \ref{cusp and swallowtail}. 
For example, the property of the equation ${\mathrm{Hess}} = 1$ affects it. 
In fact from the equations $\omega = dx \wedge dy - dp \wedge dq = 0, \theta 
= dz - pdx - qdy = 0$, we see $\pi_1\circ f$ is immersive at $x_0$ if and only if 
$\pi_2\circ f$ is immersive at $x_0$, therefore 
$\pi_1\circ f$ is singular at $x_0$ if and only if 
$\pi_2\circ f$ is singular at $x_0$. 

\ber 
{\rm 
In \cite{KRSUY}, it is shown that the generic singularities on flat surfaces 
in the hyperbolic 3-space $H^3$ 
are cuspidal edges and swallowtails, based on the representation formula obtained in \cite{KUY}. 
Also a useful criterion on singularities 
is established in \cite{KRSUY}, so called \lq\lq KRSUY" criterion. 
To show Theorem \ref{cusp and swallowtail}, we analyze by means of power series expansions and transversality arguments, in each case  
${\rm{Hess}} = 1, {\rm{Hess}} = -1, K = 1, K = -1$ and 
${\rm{Hess}} = 0 (K = 0)$ on $\R^3$. 
To finish up the classification in each case, we apply the \lq\lq KRSUY" criterion.  

For the improper affine spheres, the global complex representation is given 
in \cite{FMM}. Moreover the notion of {\it improper affine maps} 
is introduced in \cite{Mar} in connection with special Lagrangian immersions, 
and improper affine spheres with singularities are considered. 
Furthermore the classification of singular improper affine spheres can be 
reduced to the classification of singular flat fronts in $H^3$ as shown 
in \cite{KRSUY}. This is communicated 
to the first author by M. Umehara and K. Yamada. 
Note that Theorem \ref{cusp and swallowtail} gives the classification 
of singularities not only for improper affine spheres, but also for 
their affine duals as well.  

It seems to be natural to conjecture that 
the generic singularities on surfaces of constant negative Gaussian curvature 
should be cuspidal edges and swallowtails. 
(See \cite{Mc} for the pictures of singularities appearing on surfaces of 
constant negative curvature). 
Moreover, 
by numerical experiments, it can be conjectured that 
also the generic singularities on surfaces of constant positive Gaussian curvature 
should be cuspidal edges and swallowtails(\cite{KS}). 
Theorem \ref{cusp and swallowtail} answers the conjectures affirmatively. 
}
\enr 

\ber
{\rm 
Theorems \ref{cusp and swallowtail} 
describes generic singularities for both 
$\pi_1$ and $\pi_2$-projections. 
Note that it is most natural to classify geometric solutions 
under the equivalence preserving the structure of 
double fibration, namely, by posing that the contactomorphism on $M$ should 
be taken in common for $\pi_1$ and $\pi_2$. 
However, then 
it is hopeless to expect a finite list of classification 
as in Theorem \ref{cusp and swallowtail}, 
since such classification has functional moduli in general. 
}
\enr

Now recall the fundamental notions appeared in Theorem \ref{cusp and swallowtail}.  

Let $\pi : (\R^{2n+1}, 0) \to (\R^{n+1}, 0)$ be a 
Legendrian fibration. 
Two Legendrian immersions $f, g : (\R^n, 0) \to (\R^{2n+1}, 0)$ 
are called {\it Legendre equivalent} if there exist 
a contactomorphism $\Phi : (\R^{2n+1}, 0) \to (\R^{2n+1}, 0)$, 
a diffeomorphism $\sigma : (\R^n, 0) \to (\R^n, 0)$ 
and a diffeomorphism $\varphi : (\R^{n+1}, 0) \to (\R^{n+1}, 0)$  
such that the following diagram commutes: 
$$
\begin{array}{ccccc}
(\R^n, 0) & \stackrel{f}{\longrightarrow} & (\R^{2n+1}, 0) & 
\stackrel{\pi}{\longrightarrow} & (\R^{n+1}, 0) \\
{\mbox{\scriptsize{$\sigma$}}}\downarrow & { } 
& {\mbox{\scriptsize{$\Phi$}}}\downarrow & { } 
& {\mbox{\scriptsize{$\varphi$}}}\downarrow \\
(\R^n, 0) & \stackrel{g}{\longrightarrow} & (\R^{2n+1}, 0) & 
\stackrel{\pi}{\longrightarrow} & (\R^{n+1}, 0).
\end{array}
$$

Let $\pi : (\R^5, 0) \to \R^3$ be a germ of Legendrian fibration 
with respect to the contact form $\theta = dz - pdx - qdy$. 
A Legendrian immersion $f : (\R^2, 0) \to (\R^5, 0)$ is 
called a {\it cuspidal edge} (or $A_2$ briefly) with respect to $\pi$,  
if $f$ is Legendre equivalent to 
$$
(x, y, z, p, q) = (u, v^2, \dfrac{2}{3}v^3, 0, v). 
$$
In this case we say that $\pi\circ f$ has the cuspidal edge at $0$. 
A Legendrian immersion $f : (\R^2, 0) \to (\R^5, 0)$ is 
called a {\it swallowtail} (or $A_3$ briefly) with respect to $\pi$, 
if $f$ is Legendre equivalent to 
$$
(x, y, z, p, q) = (u, v^3 + uv, \dfrac{3}{4}v^4 + \dfrac{1}{2}uv^2, - \dfrac{1}{2}v^2, v). 
$$
In this case we say that $\pi\circ f$ has the swallowtail at $0$. 

The immersion can be called \lq\lq $A_1$-singularity". 

\bee
{\rm 
Let $f : (\R^2, 0) \to \R^5$ be a map-germ defined by 
$$
\begin{array}{rcl}
x & = & u \\
y & = & 2uv + 3u^2v - v^3 \\
z  & = &  
\dfrac{1}{3}u^3 + uv^2  + \dfrac{1}{4}u^4 + \dfrac{3}{2}u^2v^2 - \dfrac{3}{4}v^4 \\
p & = & u^2 - v^2 + u^3 - 3uv^2 \\
q & = & v. 
\end{array}
$$
Then $f$ is a geometric solution to ${\mathrm{Hess}} = 1$. 
Moreover $\pi_1\circ f = (x, y, z)$ has the swallowtail at $0$ and 
$\pi_2\circ f = (x, y, px + qy - z)$ has the cuspidal edge at $0$. See Figure 1. 
Note that 
$$
\tilde{z} = px + qy - z = \dfrac{2}{3}u^3 + \dfrac{3}{4}u^4 - \dfrac{3}{2}u^2v^2 
- \dfrac{1}{4}v^4. 
$$
}
\ene
\begin{figure}[htbp]
  \begin{center}
      \includegraphics[width=5truecm, clip, keepaspectratio]{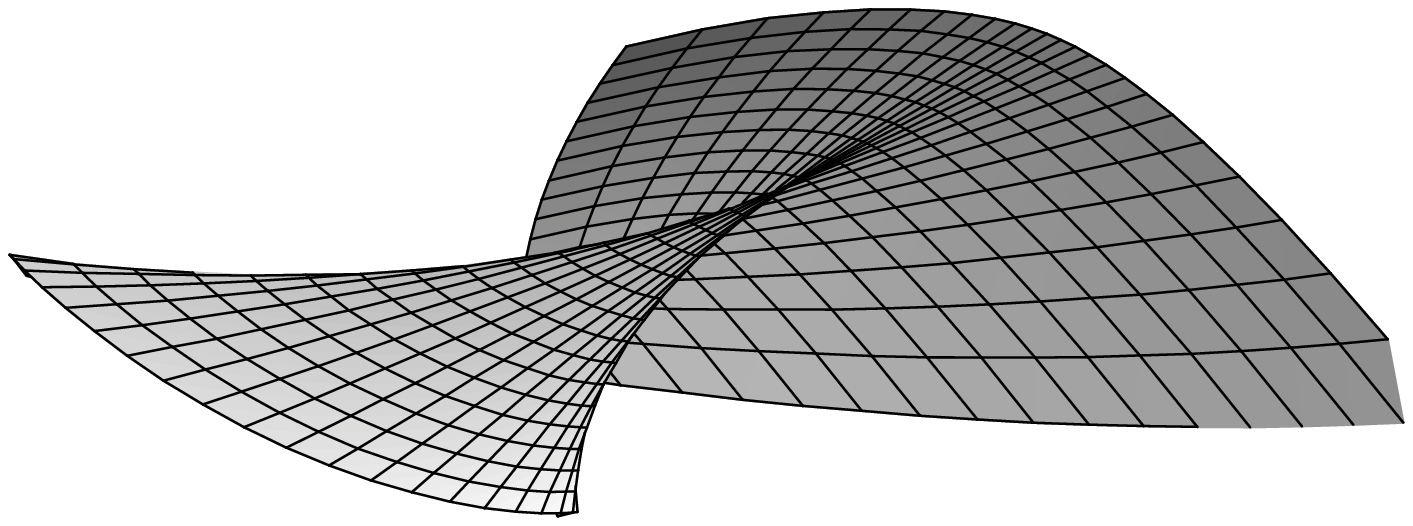} 
\quad \quad \quad \quad  
\includegraphics[width=5truecm, clip, keepaspectratio]{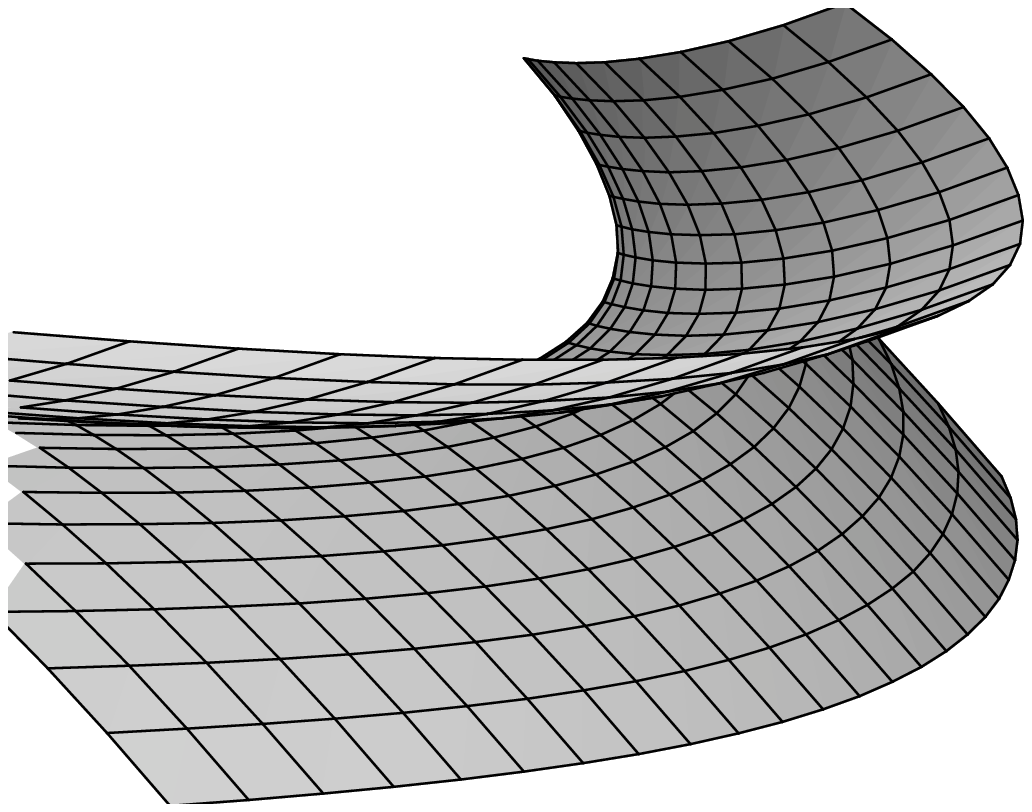} 
    \caption{An improper affine sphere with the swallowtail singularity 
(left) and with  
the cuspidal edge singularity (right) in its dual. }
    \label{integral}
  \end{center}
\end{figure}%

\ 

The solution surfaces to the equation ${\mathrm{Hess}} = 0$ or $K = 0$, in the case 
$c = 0$,  
are so called \lq\lq developable surfaces" \cite{I1}. 
Then we face a different situation with respect to dual surfaces from the cases $c \not= 0$. 
Actually we have the following classification result: 

\bet
\label{dev-front}
A generic geometric solution to ${\mathrm{Hess}} = 0$ 
{\rm (}resp. $K = 0${\rm )} has only cuspidal edges {\rm (}$A_2${\rm )} and swallowtails {\rm (}$A_3${\rm )} 
as singularities with respect to $\pi_1$, while it collapses to 
a generic immersed space curve by $\pi_2$. 
\ent 

\ 

Note that the same classification result for generic singularities of developable surfaces 
has been given 
in \cite{Izumiya}, with respect to the topology on the space of 
tangent developables to space curves. 

\ 

Also we have a result on the geometric solutions which are not necessarily immersions; 
generalized geometric solutions. 
We call a $C^\infty$ mapping $f : N^n \to M^{2n+1}$, which is not necessarily 
an immersion, an {\it integral mapping} if $f_*(TN) \subset D(\subset TM)$. 
If $D = \{ \theta = 0\}$, then the condition means that $f^*\theta = 0$. 
An integral mapping is called 
a {\it generalized geometric solution} to a Monge-Amp\`{e}re system
${\mathcal M}$ generated by a contact form $\theta$ 
and an $n$-form $\omega$ if $f^*\theta = 0, f^*\omega = 0$. 

An integral map-germ $f : (\R^2, 0) \to M^5$ is called an {\it open umbrella} 
if $f$ is contactomorphic to 
$$
(u, v) \mapsto (x, y, z, p, q) = (u, v^2, uv^3, v^3, \dfrac{3}{2}uv), 
$$
by a contactomorphism not necessarily preserving the Legendre fibrations. 
The open umbrella appears as the Legendre lifting of the tangent developable surface 
to a space curve where the curvature does not vanish and the torsion simply vanishes 
(\cite{I1}). 
Therefore the open umbrella is a generalized geometric solution to 
${\mbox{\rm{Hess}}} = 0$ {\rm (}resp. $K = 0${\rm )} (\cite{IMo}). 

Then we show the following: 

\bep
\label{OU} 
An open umbrella cannot be a generalized geometric solution to  
${\mbox{\rm{Hess}}} = c (c \not= 0)$ {\rm (}resp. 
$K = c (c \not= 0)${\rm )}. 
\enp

\

In the next section, we recall the criterion of cuspidal edges ($A_2$) 
and swallowtails ($A_3$) established in \cite{KRSUY}. 
In \S \ref{Singularities of improper affine spheres and their duals.}, 
we deal with geometric solutions to 
the equations ${\mathrm{Hess}} = \pm 1$, in connection with 
the classical solutions to the Laplace equation and to the wave equation, 
and show Theorems \ref{cusp and swallowtail} 
for them. 
In \S \ref{Singularities of pseudo-spherical surfaces and their duals.}, 
we treat geometric solutions to the equations $K = \pm 1$, 
by means of the method of Cauchy-Kovalevskaya's type 
to complete the proof 
of Theorem \ref{cusp and swallowtail}. 
In \S \ref{Singularities of developable surfaces.}, 
we give the proofs of 
Theorem \ref {dev-front} and Proposition \ref{OU}.

\section{A criterion  for cuspidal edges and swallowtails.}
\label{A criterion  for cuspidal edges and swallowtails.}

In the recent paper \cite{KRSUY}, a simple criterion is established for 
cuspidal edge singularities and swallowtails singularities of wave fronts. 
The criterion is very easy to handle, and 
we are going to use them effectively for our classification results 
in the following sections. 
We modify the criterion slightly according to the situation to which we are going to apply. 

Let $f : (\R^2, (u_0, v_0)) \to \R^5$ be a germ of Legendrian immersion, 
$$
f(u, v) = (x(u, v), y(u, v), z(u, v), p(u, v), q(u, v)), 
$$ 
for 
the standard contact form $\theta = dz - pdx - qdy$ on $\R^5$ 
with the Legendrian fibration $\pi : \R^5 \to \R^3$, 
$\pi(x, y, z, p, q) = (x, y, z)$. Set $g = \pi\circ f$. 
Suppose $(u_0, v_0)$ is a singular point, namely a non-immersive point, of $g$. 
We call the singular point $(u_0, v_0)$ {\it non-degenerate} if $\Delta = \det\left( 
\begin{array}{cc}
x_u, x_v \\
y_u, y_v 
\end{array}
\right) : (\R^2, (u_0, v_0)) \to (\R, 0)$ is submersive. 
Note that $g_u\times g_v = \Delta\cdot (-p, -q, 1)$ so that 
the singular locus of $g$ is given by $\Delta = 0$ on 
$(\R^2, (u_0, v_0))$. 

Suppose a singular point $(u_0, v_0)$ of $g$ is non-degenerate. 
Then the singular locus can be parametrized by an immersed 
curve $\gamma : (\R, 0) \to (\R^2, (u_0, v_0))$. Moreover we see 
$g$ has rank $1$ along $\gamma(t)$ near $t = 0$. 
Then the kernel field of $g_*$ is generated by a 
non-vanishing vector field $\eta : (\R, 0) \to T\R^2$ along $\gamma$ 
so that $(g_{\gamma(t)})_*(\eta(t)) = 0$. 

Then the criterion is given as follows: 

\bep
\label{KRSUY criterion}
{\rm (Proposition 1.3 of \cite{KRSUY})} 
Let $p$ is a non-degenerate singular point of $g = \pi\circ f : (\R^2, (u_0, v_0)) 
\to \R^3$ for a Legendrian immersion $f : (\R^2, (u_0, v_0)) \to \R^5$. 
Then 
$f$ is a cuspidal edge at $(u_0, v_0)$ with respect to $\pi$ if and only if 
$\det({\gamma}'(0), \eta(0)) \not= 0$. 
On the other hand, 
$f$ is a swallowtail at $(u_0, v_0)$ with respect to $\pi$ if and only if 
$\det({\gamma}'(0), \eta(0)) = 0$ 
and 
$\dfrac{\pa}{\pa t}\det({\gamma}'(t), \eta(t))\vert_{t = 0} \not= 0$. 
\enp

\bee
{\rm 
Just to make sure, let us check that our normal form 
$(x, y, z, p, q) = (u, v^2, \dfrac{2}{3}v^3, 0, v)$ of the cuspidal edge 
satisfies the criterion. 
In this case, we have $\Delta = 2v$. We can set $\gamma(t) = (t, 0)$ and 
$\eta(t) = (0, 1)$. Therefore $\det({\gamma}'(0), \eta(0)) = 1$. 

For the normal form 
$(x, y, z, p, q) = (u, v^3 + uv, \dfrac{3}{4}v^4 + \dfrac{1}{2}uv^2, - \dfrac{1}{2}v^2, v)$ 
of the swallowtail at $(u_0, v_0) = (0, 0)$, we have 
$\Delta = 3v^2 + u$. Then $\eta(t) = (0, 1)$. For any immersion 
$\gamma(t) = (u(t), v(t))$ which parametrizes $\{ \Delta = 0\}$, we have 
$3v(t)^2 + u(t) = 0$, and $u'(0) = 0$. 
Moreover, since $v(0) = 0, v'(0) \not= 0$, we have 
$\dfrac{\pa}{\pa t}\det({\gamma}'(t), \eta(t))\vert_{t = 0} = u''(0) = 6v'(0)^2 \not= 0$.

}
\ene

\section{Singularities of improper affine spheres and their duals.} 
\label{Singularities of improper affine spheres and their duals.}

We observe  the fact that 
the ellipticity of the equation ${\mathrm{Hess}} = 1$ 
implies the \lq\lq rigidity" of solutions, controlled by 
holomorphic or harmonic functions, 
while the hyperbolicity of the equation ${\mathrm{Hess}} = -1$ 
implies the \lq\lq softness" of solutions, controlled by 
$C^\infty$ functions. 

Then the key point of the proof of Theorem \ref{cusp and swallowtail} 
lies on the fact that even in elliptic cases 
there exist enough solutions 
implying the validity of a transversality theorem. 

Before starting the detailed analysis, 
we remark that our objects are homogeneous: 
Consider the group $G$ of equi-affine transformations 
preserving the vector field $\dfrac{\pa}{\pa z}$. Then $G$ acts transitively 
on $\R^5$, on $\R^3$, and the dual $\R^3$ in the natural way. 
Moreover $\pi_1 : \R^5 \to \R^3$ and 
$\pi_2 : \R^5 \to \R^3$ 
are $G$-equivariant. Furthermore 
the differential system $\omega = 0, \theta = 0$ on $\R^5$ 
is $G$-invariant. 
Throughout this section, we use the homogeneity for simplifying 
the calculations. 

\ 

First, let us consider the equation ${\mathrm{Hess}} = 1$ on $\R^5$. 

Let $f : (\R^2, 0) \to \R^5$ be a
germ of a geometric solution to ${\mathrm{Hess}} = 1$. 
We assume $f$ is an immersion, and
$f^*\theta = 0$ for the contact form
$\theta = dz - pdx - qdy$ and
$f^*\omega = 0$ 
for the $2$-form 
$\omega = dx \wedge dy - dp\wedge dq$. 
Set 
$$
f(u, v) = (x(u, v), y(u, v), z(u, v), p(u, v), q(u, v)).
$$

First we observe the following: 

\bel
\label{sing 1}
Suppose either $(\pi_1\circ f)_* : T_0\R^2 \to T_{0}\R^3$ or
$(\pi_2\circ f)_* : T_0\R^2 \to T_{0}\R^3$ is injective.
Then both $(\pi_1\circ f)_*$ and $(\pi_2\circ f)_*$ are
injective. 
\enl

\Proof
Since $f^*(dz - pdx - qdy) = 0$, 
the condition that 
$(\pi_1\circ f)_* : T_0\R^2 \to T_0\R^3$ 
is injective is equivalent to that 
$(\pi_1\circ f)^*(dx \wedge dy) \not= 0$ 
at $0$. 
Similarly, since 
$f^*(d(xp + yq - z) - xdp - ydq) = 0$, 
the condition that 
$(\pi_2\circ f)_* : T_0\R^2 \to T_0\R^3$ 
is injective is equivalent to that 
$(\pi_2\circ f)^*(dp \wedge dq) \not= 0$. 
On the other hand, $0 = f^*(dx \wedge dy - dp \wedge dq) 
=  (\pi_1\circ f)^*(dx \wedge dy) - (\pi_2\circ f)^*(dp \wedge dq)$, 
we see 
$(\pi_1\circ f)^*(dx \wedge dy) = (\pi_2\circ f)^*(dp \wedge dq)$. 
Thus we have the result. 
\QED

\ 

Moreover we have: 

\bel
\label{rank 1}
Suppose both $(\pi_1\circ f)_*$ and $(\pi_2\circ f)_*$
are not injective at $0$. Then
$(\pi_1\circ f)_*$ has rank $1$ and $(\pi_2\circ f)_*$ has rank $1$ 
at $0$.
\enl

\Proof
Assume that $(\pi_1\circ f)_*$ is not injective and dose not have rank $1$. 
Then $(\pi_1\circ f)_*$ has rank $0$, which 
means that $f_*(T_0\R^2)$ is contained in $E_1$. 
Since $(E_1)_{f(0)}$ projects to $T_{(\pi_2\circ f)(0)}\R^3$ 
injectively by $(\pi_2)_*$, we see $(\pi_2\circ f)_*$ 
must be injective. 
Thus we see $(\pi_1\circ f)_*$ has rank $1$. 
By the symmetric argument, we have also that 
$(\pi_2\circ f)_*$ has rank $1$. 
\QED

\ 

Furthermore we have: 

\bel
\label{rank 2}
Suppose both $(\pi_1\circ f)_*$ and $(\pi_2\circ f)_*$
are not injective at $0$. Then we have 
$(x\circ f, q\circ f)_* : T_0\R^2 \to T_0\R^2$ is isomorphic or 
$(y\circ f, p\circ f)_* : T_0\R^2 \to T_0\R^2$ is isomorphic. 
\enl

\Proof 
By Lemma \ref{rank 1}, we see $d(x\circ f)(0) \not= 0$ 
or $d(y\circ f)(0) \not= 0$. 
Let $d(x\circ f)(0) \not= 0$. 
Since $f$ is a Legendrian immersion, we see 
$(x, p), (x, q), (y, p)$ or $(y, q)$ composed with $f$ 
is a local diffeomorphism. 
Assume $d(x\circ f)$ and $d(p\circ f)$ 
are linearly independent at $0$. 
Since $f^*(d\theta) = 0$, we see 
$0 \not= d(p\circ f) \wedge d(x\circ f) = - d(q\circ f) \wedge 
d(y\circ f)$ at $0$. 
Therefore $d(q\circ f)(0)$ and $d(y\circ f)(0)$ 
must be linearly independent. 
Besides, we are supposing that $d(x\circ f)(0)$ 
and $d(y\circ f)(0)$ are linearly dependent. 
Therefore we have that 
$d(x\circ f)(0)$ and $d(q\circ f)(0)$ are linearly independent. 
If $d(y\circ f)(0) \not= 0$, then similarly we have 
$d(y\circ f)(0)$ and $d(p\circ f)(0)$ are linearly independent. 
\QED

Now suppose 
either $(\pi_1\circ f)_* : T_0\R^2 \to T_0\R^3$ or
$(\pi_2\circ f)_* : T_0\R^2 \to T_0\R^3$ is not injective.
By Lemmata \ref{sing 1}, \ref{rank 1}, \ref{rank 2}, 
we may assume  $x = u, q = v$ 
by a coordinate transformation of $(\R^2, 0)$ 
and by the isomorphism $(x, y, z, p, q) \mapsto (y, x, z, q, p)$ 
of the system. Thus we set 
$$
f(u, v) = (u, y(u, v), z(u, v), p(u, v), v).
$$
By the condition $dz = pdu + vdy$, we have
that $dp \wedge du + dv \wedge dy = 0$.
Then we have 
$$
\dfrac{\pa p}{\pa v} = - \dfrac{\pa y}{\pa u}.
$$
Also by the condition $\omega = dp \wedge dv - du \wedge dy = 0$,
we have 
$$
\dfrac{\pa p}{\pa u} = \dfrac{\pa y}{\pa v}.
$$
Thus we obtain the Cauchy-Riemann equation. Therefore
we see the complex valued function
$p + \sqrt{-1} y : (\R^2, 0) \to \C$
is holomorphic with the complex coordinate $u + \sqrt{-1}v$.

\ber
\label{partialL}
{\rm 
Consider a contact transformation $L : \R^5 \to \R^5$ 
defined by $L(x, y, z, p, q) 
= (x, q, z - yq, p, -y)$, 
which is called a {\it partial Legendre transformation}. 
Then $L^*(dp \wedge dy - dq \wedge dx) 
= - (dx \wedge dy - dp\wedge dq) = -\omega$. 
Therefore, as well known, the equation ${\mathrm{Hess}} = 1$ 
is contactomorphic to the Laplace equation $z_{xx} + z_{yy} = 0$. 
Similarly, we see the equation ${\mathrm{Hess}} = -1$ 
is contactomorphic to the wave equation $z_{xx} - z_{yy} = 0$. 
Thus the \lq\lq nonlinear" equation ${\mathrm{Hess}} = \pm 1$, or 
$rt - s^2 = \pm 1$ in the classical Goursat's notation, is 
equivalent to the \lq\lq linear" equation $r \pm t = 0$. 
Actually we have used this procedure of Legendrian transformation, 
and naturally we have got the Cauchy-Riemann 
equation (resp. the wave equation). 
However note that the Legendrian transformation 
$L$ does not preserve the structure of the Lagrangian pair, 
and also that the Cauchy-Riemann equation we have got is 
not over the $xy$-plane but rather over the 
$uv$-plane, the parameter plane 
of the geometric solution. 
}
\enr

Thus we have: 

\bep
Let $f : (\R^2, 0) \to \R^5$ be
a germ of a geometric solution of the 
Monge-Amp\`{e}re system 
associated to the equation ${\mathrm{Hess}} = 1$.
Then there exists a germ of a holomorphic function 
$h = p + \sqrt{-1}y : (\R^2, 0) = (\C, 0) \to (\C, 0)$ of $u + \sqrt{-1}v$ 
such that 
$f$ is Legendre equivalent to
$f_h = (u, y(u, v), z(u, v), p(u, v), v)$, 
up to a diffeomorphism on $(\R^2, 0)$ 
and a contactomorphism on $\R^5$ preserving $({\mathcal{M}}, E_1, E_2)$, 
defined by the line integral 
$$
z (u, v) = \int_{(0,0)}^{(u,v)} \left( p(u, v) + v\dfrac{\pa y}{\pa
u}\right) 
du + v\dfrac{\pa y}{\pa v}dv.
$$
\enp

\bee
{\rm 
Consider the holomorphic function 
$$
h  =  (u + \sqrt{-1}v)^2  =  u^2 - v^2 + \sqrt{-1}(2uv). 
$$ 
Set $p = u^2 - v^2$ and
$y = 2uv$. Then we have 
$
f_h = (u, 2uv, \dfrac{1}{3}u^3 + uv^2, u^2 - v^2, v).
$
For the function $h = \sqrt{-1}(u + \sqrt{-1}v)^2$ we have
$
f_h = (u, u^2 - v^2, -\dfrac{2}{3}v^3, -2uv, v).
$

More generally, for a holomorphic function 
$$
h =  (a_1 + \sqrt{-1}b_1)(u + \sqrt{-1}v) 
+ (a_2 + \sqrt{-1}b_2)(u + \sqrt{-1}v)^2 
+ (a_3 + \sqrt{-1}b_3)(u + \sqrt{-1}v)^3 + \cdots, 
$$
we have $f_h = (x, y, z, p, q)$ with 
$$
\begin{array}{rcl}
x & = & u \\
y & = &  a_1v + b_1u + a_2(2uv) + b_2(u^2 - v^2) + a_3(3u^2v - v^3) + b_3(u^3 - 3uv^2) + \cdots \\
z & = & a_1(\dfrac{1}{2}u^2 + \dfrac{1}{2}v^2) + 
a_2(\dfrac{1}{3}u^3 + uv^2) + b_2(-\dfrac{2}{3}v^3) \\
{ } & { } & 
+ a_3(\dfrac{1}{4}u^4 + \dfrac{3}{2}u^2v^2 - \dfrac{3}{4}v^4) + b_3(-2uv^3) + \cdots \\
p & = & a_1u - b_1v + 
a_2(u^2 - v^2) + b_2(-2uv) + a_3(u^3 - 3uv^2) + b_3(-3u^2v + v^3) + \cdots \\
q & = & v. 
\end{array}
$$
}
\ene

This example parametrizes all geometric solutions to ${\mathrm{Hess}} = 1$ 
with projection singularities, 
up to the equivalence and up to their $3$-jets. 
Namely, for any geometric solution, its germ at each point can be represented as above 
for some coordinates $u, v$ centered at the point. 

To be precise, the argument goes as follows: 
Let $U$ be an open subset in the $uv$-plane $\R^2$. 
We identify the $3$-jet space $J^3(U, \R^3)$ with 
the submanifold of $J^3(U, \R^5)$ consisting of 
$3$-jets $j^3f(x_0)$ with a form 
$$
f(u, v) = (u, y(u, v), z(u, v), p(u, v), v).
$$
Moreover we identify $J^3(2, 3)$ with the submanifold of 
$J^3(2, 5) = \{ j^3g(0) \mid g : (\R^2, 0) \to (\R^5, 0) \}$ 
consisting of $3$-jets  $j^3g(0)$ with $x\circ g = u$ and 
$q\circ g = v$. 
Consider the submersion 
$$
\Phi : J^3(U, \R^5) \to J^3(2, 5)
$$
by $\Phi(j^3f(u_0, v_0)) = j^3(\bar{f}(0))$, 
where $\bar{f} : (\R^3, 0) \to (\R^5, 0)$ is defined by 
$\bar{f}(u, v) = f(u+u_0, v+v_0) - f(u_0, v_0)$. 
Then $\Phi$ maps $J^3(U, \R^3)$ onto $J^3(2, 3)$. So, the submersion 
$$
\Phi : J^3(U, \R^3) \to J^3(2, 3)
$$
is induced. 
Let ${\mathcal G} = {\mathcal G}_U$ be 
the $3$-jet space of geometric solutions $(U, (u_0, v_0)) \to \R^5$ with a form 
$$
f(u, v) = (u, y(u, v), z(u, v), p(u, v), v).
$$
Moreover we identify $\R^6$ with the submanifold in $J^3(2, 3)$ consisting of $j^3(f_h)(0)$ described as above.  
Then we have seen 
$\Phi^{-1}(\R^6) = {\mathcal G}$. 
Thus we see ${\mathcal G}$ is a manifold of dimension $11$ and 
$\Phi : {\mathcal G} \to \R^6$ is a submersion. 

Now the proof of Theorem \ref{cusp and swallowtail} for the equation 
${\mbox{\rm{Hess}}} = 1$ is achieved as follows: 
Consider the family of $3$-jets of geometric solutions 
of the above form parametrized by parameters $a_1, a_2, a_3, b_1, b_2, b_3$. 

The singular locus of $\pi_1\circ f$ is described by 
$$
\Delta = y_v = a_1 + 2a_2u - 2b_2v + 3a_3(u^2 - v^2) - 6b_3uv + \cdots = 0, 
$$
while also the singular locus of $\pi_2\circ f$ is described by 
$$
\Delta = p_u = p_v = a_1 + 2a_2u - 2b_2v + 3a_3(u^2 - v^2) - 6b_3uv + \cdots = 0. 
$$

Consider the hypersurface $\{ a_1 = 0 \}$ in $\R^6$. Then 
$\Phi^{-1}(\{ a_1 = 0 \}) \subset {\mathcal G}$ is also a hypersurface. 

Let $f_h : U \to \R^5$ be the geometric solution to 
${\mbox{\rm{Hess}}} = 1$ defined by a 
holomorphic function $h : U (\subset \R^2 = \C) \to \C$. 
Consider the mapping $\Psi(h) : U \to {\mathcal G}$ 
defined by 
$\Psi(h)(u_0, v_0) = j^3f_h(u_0, v_0)$. 
Then,  by a small perturbation 
$\tilde{h}(z) = h(z) + \alpha z + \beta z^2 + \gamma z^3$ by a complex polynomial 
$(\alpha, \beta, \gamma \in \C)$, 
we can make $\Psi(\tilde{h}) : U \to {\mathcal G}$ transverse to 
$\Phi^{-1}(\{ a_1 = 0 \})$. Then 
the locus $\{ a_1 = 0 \}$ is a smooth curve in $U$. 

If (i) $a_1\not= 0$, then both $\pi_1\circ f$ and $\pi_2\circ f$ are immersive. 

Moreover by a perturbation of $h$ if necessary, 
we may suppose, along the curve $\{ a_1 = 0\}$, 
there occurs only three cases: 
(ii) $a_1 = 0, a_2 \not= 0, b_2 \not= 0$, 
(iii) $a_1 = 0, a_2 \not= 0, b_2 = 0, a_3 \not= 0, $ or 
(iv) $a_1 = 0, a_2 = 0, b_2 \not= 0, a_3 \not= 0. $

We apply the criterion Proposition \ref{KRSUY criterion} to our case 
to verify the 
singularity type of $\pi_1\circ f$ (resp. $\pi_2\circ f$). 

The singularities are non-degenerate since $(a_2, b_2) \not= (0, 0)$. 

Now suppose (iii) $a_1 = 0, a_2 \not= 0, b_2 = 0, a_3 \not= 0$. 
Then we have 
$$
\begin{array}{rcl}
x & = & u, \\
y & = & b_1u + a_2(2uv) + a_3(3u^2v - v^3) + b_3(u^3 - 3uv^2) + \cdots, \\
z  & = &  
a_2(\dfrac{1}{3}u^3 + uv^2)  + a_3(\dfrac{1}{4}u^4 + \dfrac{3}{2}u^2v^2 - \dfrac{3}{4}v^4) + b_3(-2uv^3) + \cdots \\
p & = & -b_1v + a_2(u^2 - v^2) + a_3(u^3 - 3uv^2) + b_3(-3u^2v + v^3) + \cdots \\
q & = & v. 
\end{array}
$$

Then we have, for $(u_0, v_0) = (0, 0)$, 
$$
y_v = 2a_2u + 3a_3(u^2 - v^2) - 6b_3uv + \cdots. 
$$
Take a parametrization $\gamma(t) = (u(t), v(t))$ of the singular locus 
$\{ y_v = 0\}$ of $\pi_1\circ f$ satisfying $\gamma(0) = (0, 0), 
\gamma'(0) \not= (0, 0)$.  
Then we see $u'(0) = 0$
$u''(0) \not= 0$ if and only if $a_3 \not= 0$. 
Since we can take $\eta(t) = (0, 1)$ as a kernel field, 
we have $\det(\gamma'(0), \eta(0)) = u'(0) = 0$ and 
$\dfrac{\pa}{\pa t}\det(\gamma'(t), \eta(t))\vert_{t = 0} = u''(0) \not= 0$. 
On the other hand, for a parametrization $\gamma(t) = (u(t), v(t))$ 
of the singular locus $\{ p_u = 0\}$ of $\pi_2\circ f$, 
$$
p_u = 2a_2u + 3a_3(u^2 - v^2) - 6b_3uv + \cdots, 
$$
we have $u'(0) = 0$, so $v'(0) \not= 0$. 
Since we can take $\eta(t) = (1, 0)$ as a kernel field for $\pi_2\circ f$, 
we have $\det(\gamma'(0), \eta(0)) = v'(0) \not= 0$

Therefore, by Proposition \ref{KRSUY criterion}, 
$\pi_1\circ f$ is the swallowtail 
and $\pi_2\circ f$ is the cuspidal edge in the case (iii). 
 
In the case (iv) $a_1 = 0, b_1 = 0, a_2 = 0, b_2 \not= 0, a_3 \not= 0, $ 
we see, similarly to the case (iii), $\pi_1\circ f$ is the cuspidal edge 
and $\pi_2\circ f$ is the swallowtail. 
Moreover in the case (ii) $a_2 \not= 0, b_2 \not= 0$, we see 
both $\pi_1\circ f$ and $\pi_2\circ f$ are the cuspidal edges. 

Thus the proof of Theorem \ref{cusp and swallowtail} for the equation 
${\mbox{\rm{Hess}}} = 1$ is completed. 

\ 

Next we consider the equation ${\mathrm{Hess}} = -1$. 
In this case we set 
$$
\theta = dz - pdx - qdy,  \quad \omega = dx \wedge dy + dp \wedge dq. 
$$
Let $f : (\R^2, 0) \to \R^5$ be a germ of an geometric solution 
to ${\mathcal{M}} = \langle \theta, \omega\rangle$. 
Then up to equivalence, we can write 
$$
f(u, v) = (u, y(u, v), z(u, v), p(u, v), v), 
$$
and we have the equation: 
$$
\dfrac{\pa p}{\pa v} = - \dfrac{\pa y}{\pa u}, \quad 
\dfrac{\pa p}{\pa u} = - \dfrac{\pa y}{\pa v}.
$$
For this wave equation, 
we get 
$$
y = \varphi(u+v) + \psi(u-v), \quad 
p = - \varphi(u+v) + \psi(u-v), 
$$
for smooth functions $\varphi, \psi$. 

If 
$$
\begin{array}{rcl}
\varphi & = & 
\varphi_1 t + \varphi_2 t^2 + \varphi_3 t^3 + \cdots, \\
\psi & = & 
\psi_1 t + \psi_2 t^2 + \psi_3 t^3 + \cdots, 
\end{array}
$$
where $\varphi_i, \psi_j$ are real numbers, 
then we have the expansions 
$$
\begin{array}{rcl}
y & = & 
\varphi_1(u+v) + \psi_1(u-v) 
+ \varphi_2(u+v)^2 + \psi_2(u-v)^2 
+ \varphi_3(u+v)^3 + \varphi_3(u-v)^3 + \cdots \\ 
p & = & 
- \varphi_1(u+v) + \psi_1(u-v) 
- \varphi_2(u+v)^2 + \psi_2(u-v)^2 
- \varphi_3(u+v)^3 + \varphi_3(u-v)^3 + \cdots, 
\end{array}
$$
up to $3$-jets. 
Moreover  we have 
$$
\begin{array}{rcl}
z & = & {\displaystyle{\int_{(0, 0)}^{(u, v)}}} pdx + qdy \\ 
{ } & =  & \dfrac{1}{2}\varphi_1(-u^2+v^2) + \dfrac{1}{2}\psi_1(u^2-v^2) 
+ \dfrac{1}{3}\varphi_2(-u^3+3uv^2+2v^3) + \dfrac{1}{3}\psi_2(u^3-3uv^2-2v^3) \\
{ } & { } & 
+ \dfrac{1}{4}\varphi_3(-u^4+6u^2v^2+8uv^3+3v^3) 
+ \dfrac{1}{4}\psi_3(u^4-6u^2v^2+8uv^3-3v^3) + \cdots. 
\end{array}
$$
Then we see Theorem \ref{cusp and swallowtail} holds 
for the Monge-Amp\`{e}re system of 
${\mathrm{Hess}} = -1$ by the same way.

\ber
{\rm 
To describe the jets of solutions to ${\mathrm{Hess}} = -1$, 
we can also expand as a \lq\lq split holomorphic function"
$$
p - jy = (a_1 + jb_1)(u + jv) + (a_2 + jb_2)(u + jv)^2 + 
(a_3 + jb_3)(u + jv)^3 + \cdots, 
$$
using the {\it split complex number} 
$u + jv$, where $j^2 = 1$, not $j^2 = -1$, and $u, v$ are real numbers, 
so as 
$$
\begin{array}{rcl}
p & = & a_1u + b_1v + a_2(u^2 + v^2) + b_2(2uv) + 
a_3(u^3 + 3uv^2) + b_3(3u^2v + v^3) + \cdots,  \\
y & = & - a_1v - b_1u - a_2(2uv) - b_2(u^2 + v^2) - a_3(3u^2v + v^3) 
- b_3(u^3 + 3uv^2) + \cdots. 
\end{array}. 
$$
In fact the coefficients are related by 
$a_k = - \varphi_k + \psi_k, b_k = - \varphi_k - \psi_k$. 
Nevertheless we have to remark that the equation ${\mathrm{Hess}} = -1$ 
admits infinitely flat nonzero perturbations of geometric solutions, while 
\lq\lq the theorem of identity" holds for solutions to ${\mathrm{Hess}} = 1$. 
}
\enr

\section{Singularities of surfaces of constant Gaussian curvature and their duals.} 
\label{Singularities of pseudo-spherical surfaces and their duals.}

Now we turn to study on geometric solutions to the equation $K = c \, (c \not= 0)$ in 
$\R^3\times S^2$.

Recall that an immersion $f : N \to \R^3\times S^2$, 
$$
f(u, v) = (x_1(u, v), x_2(u, v), x_3(u, v), y_1(u, v), y_2(u, v), y_3(u, v)), 
$$
is called a geometric solution to $K = c$ if $f$ satisfies 
the conditions $f^*\theta = 0, f^*\omega = 0$ and of course 
$y_1^2 + y_2^2 + y_3^2 = 1$, 
where 
$$
\begin{array}{rcl}
\theta & = & y_1dx_1 + y_2dx_2 + y_3dx_3, \\
\omega & = & 
c(y_1dx_2\wedge dx_3 + y_2dx_3\wedge dx_1 + y_3dx_1\wedge dx_2) \\
{ } & { } & \qquad 
- (y_1dy_2\wedge dy_3 + y_2dy_3\wedge dy_1 + y_3dy_1\wedge dy_2). 
\end{array}
$$

First we remark that the Euclidean 
group $G$ on the Euclidean space $\R^3$ 
acts also on the unit tangent bundle $\R^3\times S^2$ 
and $\R\times S^2$ transitively such that 
$\pi_1 : \R^3\times S^2 \to \R^3$ and 
$\pi_2 : \R^3\times S^2 \to \R\times S^2$ 
are both $G$-equivariant. 
Moreover the Monge-Amp\`ere system associated to 
the equation $K = c$ is also $G$-invariant. 
For each $(x_0, y_0) \in \R^3\times S^2$, 
the quotient mapping 
$\pi : G \to \R^3\times S^2$, $\pi(g) = g\cdot(x_0, y_0)$, 
is a $C^\infty$ fibration. 
Then there exists a local $C^\infty$ section $S = S(x_0, y_0; \cdot, \cdot ) 
: V \to G$ over a neighborhood $V$ of 
$(x_0, y_0)$. 
Note that $S(x_0, y_0, x, y) \in G$ transforms 
$(x_0, y_0)$ to $(x, y)$ for any $(x, y) \in V$.

Let $f : (N, x_0) \to \R^3\times S^2$ be a 
germ of a geometric solution to $K = c$. 
Take a system of coordinates $(u, v)$ 
centered at $x_0$ of $N$, and fix 
a $g_0 \in G$ transforming $f(x_0)$ to 
$b = (0, 0, 0; 1, 0, 0) \in \R^3\times S^2$. 
Then, for each $(u_0, v_0)$ near $(0, 0)$, 
we define $\bar{f}_{u_0, v_0} : (\R^2, (0, 0))
\to (\R^3\times S^2, b)$ by 
$$
\bar{f}_{u_0, v_0}(u, v) = g_0\cdot S(f(0, 0) ; f(u_0, v_0))^{-1}\cdot f(u + u_0, 
v + v_0). 
$$

Thus, by the homogeneity of the equation $K = c$, 
we may suppose that $f(x_0) = b = (0, 0, 0; 1, 0, 0)$.  

Set $p = -y_2/y_1, q = - y_3/y_1$. Then we have 
$$
\omega = y_1\{ c(dx_2\wedge dx_3 - pdx_3
\wedge dx_1 - qdx_1\wedge dx_2) - y_1^2 dp\wedge dq\}, 
$$
and $y_1^2 = \dfrac{1}{1+p^2+q^2}$. 
Then the Monge-Amp\`{e}re system for $K = c$ is locally given by 
$$
\left\{ 
\begin{array}{l}
c(dx_2\wedge dx_3 - pdx_3\wedge dx_1 - qdx_1\wedge dx_2) - \dfrac{1}{1+p^2+q^2}dp\wedge dq = 0, \\
dx_1 - pdx_2 - qdx_3 = 0. 
\end{array}
\right.
$$
Setting $x_1 = z, x_2 = x, x_3 = y$, we have 
$$
\left\{ 
\begin{array}{l}
c(1 + p^2 + q^2)^2 dx\wedge dy - dp\wedge dq = 0, \\
dz - pdx - qdy = 0. 
\end{array}
\right.
$$
Then, for some local coordinates on $\R^3$ and $\R\times S^2$, 
$\pi_1$ and $\pi_2$ are given by 
$\pi_1(x, y, z, p, q) = (x, y, z)$ and $\pi_2(x, y, z, p, q) = (xp + yq - z, p, q)$ 
respectively. 

As in the case Hess $= c$,  we may suppose the mapping 
$(u, v) \mapsto (x(u, v), q(u, v))$ is a local diffeomorphism. So 
we assume $x = u$ and $q = v$. 

Now, we consider the partial Legendre transformation 
$$
L(x, y, z, p, q) = (x, q, z - yq, p, -y), 
$$
and its inverse 
$$
L^{-1}(x, y, z, p, q) = (x, -q, z - yq, p, y)
$$
(cf. Remark \ref{partialL}). 
Then $L\circ f : (\R^2, 0) \to (\R^5, 0)$, 
$$
L\circ f(u, v) = (u, v, Z, P, Q) 
$$
satisfies the equation 
$$
\left\{ 
\begin{array}{l}
c(1 + P^2 + v^2)^2 dQ\wedge du - dP\wedge dv = 0, \\
dZ - Pdu - Qdv = 0. 
\end{array}
\right.
$$
Note that 
$$
Z(u, v) = z\circ f(u, v)  - (y\circ f(u, v))v, \ 
P(u, v) = p\circ f(u, v), \ 
Q(u, v) = - y\circ f(u, v). 
$$
Then we have 
$$
c(1 + P^2 + v^2)^2 Q_v + P_u = 0, \ 
P = Z_u, \ Q = Z_v, 
$$ 
so we have the Monge-Amp\`{e}re equation 
$$
Z_{uu} + c(1 + Z_u^2 + v^2)^2 Z_{vv} = 0 \qquad (*), 
$$
on a function $Z = Z(u, v)$ 
with $Z(0, 0) = 0, Z_u(0, 0) = 0, Z_v(0, 0) = 0$. 
Thus we reduce the problem on geometrical solutions (with projection singularities), 
via a Legendrian transformation, 
to the problem on classical solutions to another Monge-Amp\`{e}re equation. 

We compute the \lq\lq prolongations" of the equation (*) to 
obtain the Taylor expansion of $Z$. 
We have, by setting $(u, v) = (0, 0)$, 
$$
Z_{uu}(0, 0) + c \, Z_{vv}(0, 0) = 0. 
$$
By differentiating by $u$ (resp. by $v$) of both sides of (*) and by 
setting $(u, v) = (0, 0)$, we have 
$$
Z_{uuu}(0, 0) + c \,  Z_{uvv}(0, 0)  = 0,  \quad 
Z_{uuv}(0, 0) + c \,  Z_{vvv}(0, 0)  = 0. 
$$
By differentiating by $u$ or $v$ of both sides of (*) twice and by 
setting $(u, v) = (0, 0)$, we have 
$$
\begin{array}{rcl}
Z_{uuuu}(0, 0) + 4c \,  Z_{uu}(0, 0)^2 Z_{vv}(0, 0) 
+ c \, Z_{uuvv}(0, 0) & = & 0, \\
Z_{uuuv}(0, 0) + 4c \,  Z_{uu}(0, 0)Z_{uv}(0, 0)Z_{vv}(0, 0)  
+ c \, Z_{uvvv}(0, 0) & = & 0, \\
Z_{uuvv}(0, 0) + 4c (Z_{uv}(0, 0)^2 + 1)Z_{vv}(0, 0)  
+ c \, Z_{vvvv}(0, 0) & = & 0. 
\end{array}
$$
If we set  
$$
\begin{array}{rcl}
Z(u, v)  & = & \dfrac{1}{2}Au^2 + Buv + \dfrac{1}{2}Cv^2 + 
\dfrac{1}{6}Du^3 + \dfrac{1}{2}Eu^2v + \dfrac{1}{2}Fuv^2 + \dfrac{1}{6}Gv^3  \\
& & \qquad + \dfrac{1}{24}Hu^4 + \dfrac{1}{6}Iu^3v + \dfrac{1}{4}Ju^2v^2 
+ \dfrac{1}{6}Kv^3 + \dfrac{1}{24}Lv^4 + \cdots, 
\end{array}
$$
then we have 
$$
\begin{array}{c}
A + c\, C = 0, \ D + c\, F = 0, \ E + c\, G = 0, \\ 
H + 4c\, A^2C + c\, J = 0, \ 
I + 4c\,  ABC + c\, K = 0, \ 
J + 4c(B^2 + 1)C + c\, L = 0. 
\end{array}
$$
Therefore we have 
$$
\begin{array}{c}
A = - c\, C, \ D = - c\, F, \ E = -c\, G, \\ 
I = 4c^2 BC^2 - c\, K, \ J = -4c(B^2 + 1)C - c\, L, 
\ H = - 4c^3 C^3 + 4c^2(B^2 + 1)C  + c^2 L.  
\end{array}
$$
Thus we have 
$$
\begin{array}{rcl}
Z(u, v)  & = & -\dfrac{c}{2}Cu^2 + Buv + \dfrac{1}{2}Cv^2 - 
\dfrac{c}{6}Fu^3 - \dfrac{c}{2}Gu^2v + \dfrac{1}{2}Fuv^2 + \dfrac{1}{6}Gv^3  \\
& & \quad + \dfrac{c^2}{24}(- 4c\, C^3 + 4(B^2 + 1)C  + L)u^4 + 
\dfrac{c}{6}(4c\, BC^2 - K)u^3v \\ 
& & \quad \quad - \dfrac{c}{4}(4B^2C + 4C + L)u^2v^2 
+ \dfrac{1}{6}Kv^3 + \dfrac{1}{24}Lv^4 + \cdots. 
\end{array}
$$
Note that the original $f$ is given by 
$$
x = u, \ y = Z_v, \ z = Z - {Z_v}v, \ p = Z_u, \ q = v. 
$$

The above procedure gives us all formal solutions to the equation (*). 
In fact, the famous {\it theorem of Cauchy-Kovalevskaya} (\cite{Jo}) 
says that the formal solution 
$Z(u, v)$ is uniquely determined when the initial conditions 
$Z(0, v)$ and $Z_u(0, v)$ are given. 
In fact, the coefficients of the Taylor expansion of $Z$ up to degree  
$r$ is determined as explicit polynomials by those of $Z(0, v)$ and $Z_u(0, v)$ 
up to degree $r$. 
Moreover, if the initial data $Z(0, v)$ and $Z_u(0, v)$ are an analytic function, then 
the solution $Z$ should be an analytic function. 

\ 

Let $Z(u, v)$ be a $C^\infty$ solution to (*) corresponding to 
a germ of a geometric solution $f$ to $K = c$. 
By taking Taylor polynomials of $Z(0, v)$ and $Z_u(0, v)$ of arbitrarily high degree, 
as initial conditions, 
we get an approximation $\tilde{Z}(u, v)$ of $Z(u, v)$ 
in $C^\infty$ topology, which is an analytic solution to (*). 

By considering $\bar{f}_{u_0, v_0} : (\R^2, 0) \to (\R^3\times S^2, b)$ 
defined as above and by considering the Taylor expansion of 
the $Z$-component of $L \circ \bar{f}_{u_0, v_0}$, we get a 
$C^\infty$ map-germ $\Phi(f) : (N, x_0) \to \R^6$, 
$$
(u_0, v_0) \mapsto (B(u_0, v_0), C(u_0, v_0), F(u_0, v_0), G(u_0, v_0), K(u_0, v_0), L(u_0, v_0)). 
$$

Because we can control the coefficients $B, C, F, G, K, L$ freely in the approximation precess, 
we can perturb $f$ into $\tilde{f}$ such that  
$\Phi(\tilde{f})$ is transverse to a given stratification (or a locally finite family of submanifolds) of $\R^6$. 

Now consider a stratification 
of $\R^6$ with coordinates 
$(B, C, F, G, K, L)$: 
$$
W^0 = \{ C \not= 0 \}, \ 
W^1 = \{ C = 0, F \not= 0, G \not= 0 \}, \ 
W_1^2 = \{ C = 0, F\not = 0, G = 0, L \not= 0 \}, 
$$
$$
W_2^2 = \{ C = 0, F = 0, G \not= 0, L \not= 0 \}, \ 
W_1^3 = \{ C = 0, F = 0, G = 0, L \not= 0 \}, \ 
$$
$$
W_2^3 = \{ C = 0, F = 0, G \not= 0, L = 0 \}, \ 
W_3^3 = \{ C = 0, F \not= 0, G = 0, L = 0 \}, 
$$
$$ 
W^4 = \{ C = 0, F = 0, G = 0, L = 0 \}. 
$$
Note that $W_1^3, W_2^3, W_3^3$ are of codimension $3$ and 
$W^4$ is of codimension $4$.

Suppose $\Phi(\tilde{f})$ is transverse to the stratification, namely, transverse 
to every $W^i_j$ in $\R^6$. 
The transversality condition 
implies that the image of $\Phi(\tilde{f})$ does not touch 
the subset $W_1^3 \cup W_2^3 \cup W_3^3 \cup W^4$, so that 
it is contained in $W^0 \cup W^1 \cup W_1^2 \cup W_2^2$. 
Then the $uv$-plane is divided into 

\begin{center}
(i) $C \not= 0$, \quad  
(ii) $C = 0, F \not= 0, G \not= 0$, \\ 
(iii) $C = 0, F = 0, G \not= 0, L \not= 0$, \quad  
(iv) $C = 0, F \not= 0, G = 0, L \not= 0$. 
\end{center}

\

Now the singular locus of $\pi_1\circ f$ is given by 
$$
Q_v = Z_{vv} = C + Fu + Gv + \dfrac{1}{2}Ju^2 
+ Kuv + \dfrac{1}{2}Lv^2 + \cdots  = 0. 
$$
We apply Proposition \ref{KRSUY criterion} to 
$\pi_1\circ f$ at $(0, 0)$. 
If $C \not= 0$, then $\pi_1\circ f$ is an immersion at $(0, 0)$. 
If $C = 0$, then $\pi_1\circ f$ is singular at $(0, 0)$. 
In the case $C = 0$, $\pi_1\circ f$ is non-degenerate at $(0, 0)$ if and only if $(F, G) \not= (0, 0)$.  
Let $\gamma(t) = (u(t), v(t))$ parametrize the singular locus 
of $\pi_1\circ f$ with $\gamma(0) = (0, 0)$. 
Then we have $Fu'(0) + Gv'(0)  = 0$. 
We set $\eta(t) = (0, 1)$. Then we have 
$\det(\gamma'(t), \eta(t)) = u'(t)$. 
Suppose $G \not= 0$. Then $u'(0) \not= 0$. 
Hence $\pi_1\circ f$ is a cuspidal edge at $(0, 0)$. 
If $C = 0, G = 0$, then we see $u'(0) = 0$. 
Moreover, if $u'(0) = 0$, then $u''(0) \not= 0$ provided 
$L \not= 0$. 
Therefore, if $C = 0, G = 0, F \not= 0, L \not= 0$, then we see 
$\pi_1\circ f$ is a swallowtail at $(0, 0)$. 

Similarly we apply Proposition \ref{KRSUY criterion} to $\pi_2\circ f$. 
The singular locus of $\pi_2\circ f$ is defined by 
$$
P_u = Z_{uu} = A + Du + Ev + \dfrac{1}{2}Hu^2 + Iuv  + 
\dfrac{1}{2}Jv^2 + \cdots = 0, 
$$
while $\eta(t) = (1, 0)$. 
Now recall that $c\not= 0$. 
Then we see that, if $C \not= 0$, then $\pi_2\circ f$ is an immersion at $(0, 0)$. 
If $C = 0, G \not= 0$, $\pi_2\circ f$ is a cuspidal edge at $(0, 0)$. 
If $C = 0, F = 0, G \not= 0, L \not= 0$, then $\pi_2\circ f$ is a swallowtail at 
$(0, 0)$. 

Therefore in the case (i), both $\pi_1\circ f$ and 
$\pi_2\circ f$ are immersions. 
In the case (ii), both $\pi_1\circ f$ and $\pi_2\circ f$ are 
cuspidal edges. 
In the case (iii), $\pi_1\circ f$ is a swallowtail and 
$\pi_2\circ f$ is a cuspidal edge. 
In the case (iv), $\pi_1\circ f$ is a cuspidal edge and 
$\pi_2\circ f$ is a swallowtail. 

Thus we have Theorem \ref{cusp and swallowtail} for $K = \pm 1$. 

\ber
{\rm 
For the case $K = c > 0$, 
the equation 
$$
Z_{uu} + c(1 + z_u^2 + v^2)^2 Z_{vv} = 0
$$
is elliptic. Therefore, by Bernstein's theorem, any solution is analytic
(\cite{B1}\cite{B2}\cite{Ho}\cite{Pe}). Hence the method of Cauchy-Kovalevskaya 
provides all geometric solutions. 

In the hyperbolic case $K = c < 0$, 
we can describe all $C^\infty$ geometric solutions to $K = c \, (c < 0)$ locally 
as follows. 

Set 
$$
\omega = c(1 + p^2 + q^2)^2 dx \wedge dy - dp \wedge dq, 
\quad 
\theta = dz - pdx - qdy. 
$$
Then we have two decompositions 
$$
\begin{array}{rcl}
- \omega + \sqrt{-c}d\theta & = & 
(\sqrt{-c}(1 + p^2 + q^2)dx + dq)\wedge (\sqrt{-c}(1 + p^2 + q^2)dy - dp), \\ 
- \omega - \sqrt{-c}d\theta & = & 
(\sqrt{-c}(1 + p^2 + q^2)dx - dq)\wedge (\sqrt{-c}(1 + p^2 + q^2)dy + dp). 
\end{array}
$$
(These decompositions are also used in the proof of Lemma \ref{never formally full}). 

Assume $f : (\R^2, 0) \to (\R^5, 0)$ is a geometric solution to 
$K = c\,  (c < 0)$. Moreover 
we assume $f^*dx(0)$ and $f^*dq(0)$ are linearly independent. 
Then we have 
$$
\begin{array}{rcl}
f^*(\sqrt{-c}(1 + p^2 + q^2)dx + dq)\wedge f^*(\sqrt{-c}(1 + p^2 + q^2)dy - dp) & = & 0, \\ 
f^*(\sqrt{-c}(1 + p^2 + q^2)dx - dq)\wedge f^*(\sqrt{-c}(1 + p^2 + q^2)dy + dp) & = & 0, 
\end{array}
$$
$f^*(\sqrt{-c}(1 + p^2 + q^2)dx + dq)(0) \not= 0$ and  
$f^*(\sqrt{-c}(1 + p^2 + q^2)dx - dq) \not= 0$. 
Thus we have two foliation ${\mathcal F}, {\mathcal F}'$ on the $uv$-plane defined by the equation $f^*(\sqrt{-c}(1 + p^2 + q^2)dx \pm dq) = 0$. 
Each leaf of ${\mathcal F}$ 
is an integral curve to 
the differential  system $D$ :  
$$ 
\left\{ 
\begin{array}{l}
\sqrt{-c}(1 + p^2 + q^2)dx + dq = 0, \\ 
\sqrt{-c}(1 + p^2 + q^2)dy - dp = 0, \\
dz - pdx - qdy = 0, 
\end{array}
\right.
$$
and each leaf of ${\mathcal F}'$ 
is an integral curve to 
the diferential system $D'$ :  
$$ 
\left\{ 
\begin{array}{l}
\sqrt{-c}(1 + p^2 + q^2)dx - dq = 0, \\ 
\sqrt{-c}(1 + p^2 + q^2)dy + dp = 0, \\
dz - pdx - qdy = 0.  
\end{array}
\right.
$$
on the $xyzpq$-space. 
Thus the geometric solution $f$ is generated by a one-parameter 
family of integral curves to $D$ (resp. $D'$). 
The differential systems $D$ and $D'$ are 
called the {\it Monge characteristic systems}. 
Compare our conclusion with  the classical  
Monge-Amp\`{e}re-Goursat theorem stating that 
any classical solution is obtained by a one-parameter family of 
integral curves to the Monge characteristic system $D$ (resp. $D'$) 
(\cite{G} Ch.2, \cite{Mat}). 
Note that the differential system 
$D \subset T\R^5$ (resp. $D' \subset T\R^5$) has 
the growth vector $(2, 3, 5)$, 
namely, $D$ is of constant rank $2$, 
$D^{(2)} = D + [D, D]$ is of constant rank $3$ and 
$D^{(3)} = D^{(2)} + [D, D^{(2)}] = T\R^5$ (cf. \cite{Mon}). 
Two foliations $({\mathcal F}, {\mathcal F}')$ on 
the $uv$-plane form 
the {\it Chebyshev net} on the image of $\pi_1\circ f$ in $\R^3$, which consists of 
{\it asymptotic lines}. Note that $\pi_1\circ f$ restricted to 
each characteristic curve is an immersion, and therefore each asymptotic line is an immersed curve 
beyond the singular locus. 
Of course, the angle $\psi$ (the {\it Chebyshev angle}) of asymptotic lines 
tends to $0$ or $\pi$ on the singular locus. 

Then $f$ is described, via a {\it B\"{a}cklund transformation}, 
by a classical solution  
$\psi = \psi(t, s)$ 
to the {\it sine-Gordon equation} $\psi_{ts} = -c \sin \psi$ (cf. \cite{PT} Ch.3). 
Here $t$ (resp. $s$) is the arclength coordinate of the characteristic curves of $D$ (resp, 
$D'$). 
Note that the system of coordinates $(t, s)$ on $\R^2$ may be different from 
the system of coordinates $(u, v)$ satisfying $x = u, q = v$. 
Then, as classically well-known, the solution $\psi$ is determined by 
two initial data $\alpha(t) = \psi(t, 0), \beta(s) = \psi(0, s)$ with 
$\alpha(0) = \beta(0)$ (\cite{G}\cite{Mat}). 
Thus we have an alternative method for the construction of 
transversal approximations of geometric solutions to $K = c\, (c < 0)$. 
For the recent progress on the initial value problem of sine-Gordon equations and 
the representation formula for surfaces with $K = c \, (c < 0)$ related to 
integrable systems, see \cite{To1}\cite{To2} for instance. 
}
\enr

\section{Singularities of developable surfaces.}
\label{Singularities of developable surfaces.}

First we show Theorem \ref {dev-front}. 
We will show that, for a generic geometric solution 
$$
f(u, v) = (x(u, v), y(u, v), z(u, v), p(u, v), q(u, v))
$$
to Hess $= 0$ (or equivalently, $K = 0$), 
$\pi\circ f$ has as singularities only cuspidal edges and swallowtails. 

First note that $\pi_2\circ f$ is of rank $\leq 1$. Thus we see 
$\pi_1\circ f$ is of rank $\geq 1$. 
Suppose $\pi_1\circ f$ is singular at a point $(u_0, v_0)$. 
Then $\pi_1\circ f$ is of rank $1$ at $(u_0, v_0)$, so that 
$\pi_2\circ f$ is of rank $1$ at $(u_0, v_0)$ as well. 
Then, up to equivalence at each such a point, we may suppose 
$$
f(u, v) = (u, y(u, v), z(u, v), p(u, v), v). 
$$
Though we can show Theorem \ref {dev-front} for $K = 0$ by the same calculations 
as in the previous section by setting $c = 0$, 
here we will give its proof in a rather direct manner. 

From $dz = pdx + qdy = pdu + vdy$ and 
$dp\wedge dq = dp \wedge dv = 0$, we have 
$$
p_v + y_u = 0, \quad p_u = 0. 
$$

From the condition $p_u = 0$, we set 
$$
p = p(v) = A + Bv + \dfrac{1}{2}Cv^2 + \dfrac{1}{6}Dv^3 + \cdots. 
$$
Then, from the condition $y_u = -p_v$,  we set 
$$
y(u, v) = 
- Bu - Cuv - \dfrac{1}{2}Duv^2 - 
\cdots 
+ \tilde{A} + \tilde{B}v + \dfrac{1}{2}\tilde{C}v^2 + \dfrac{1}{6}\tilde{D}v^3 + \cdots. 
$$
Then we have 
$$
z(u, v) = {\Tilde{\Tilde A}} + Au - \dfrac{1}{2}Cuv^2 - \dfrac{1}{3}Duv^3 - \cdots 
+ \dfrac{1}{2}\tilde{B}v^2 + \dfrac{1}{3}\tilde{C}v^3 + \dfrac{1}{6}\tilde{D}v^4 + \cdots,  
$$
and 
$$
\tilde{z}  = px + qy - z  = 
- {\Tilde{\Tilde A}} + \tilde{A}v + \dfrac{1}{2}\tilde{B}v^2 + \dfrac{1}{6}\tilde{C}v^3 + \cdots. 
$$

The singular locus of $\pi_1\circ f$ is given by 
$$
y_v =  
-  Cu - Duv - \cdots + \tilde{B} + \tilde{C}v + \dfrac{1}{2}\tilde{D}v^2 + \cdots = 0. 
$$
If $\tilde{B} \not= 0$ we see $\pi\circ f$ is immersive at $(u, v) = (0, 0)$. 
Suppose $\tilde{B} = 0$. Then $f$ is non-degenerate if and only if 
$(C, \tilde{C}) \not= (0, 0)$. 
Let $\gamma(t) = (u(t), v(t))$ parametrize the singular locus. Set 
$\eta(t) = (0, 1)$. Then $\det(\gamma'(t), \eta(t)) = u'(t)$. 
If $\tilde{B} = 0, \tilde{C} \not= 0, C \not= 0$, then $\det(\gamma'(0), \eta(0)) = 
u'(0) \not= 0$. 
Then, by Proposition \ref{KRSUY criterion}, we see 
$\pi_1\circ f$ is a cuspidal edge. 
If $\tilde{B} = 0, \tilde{C} = 0, C \not= 0$, then we have $u'(0) = 0$. 
Since $u''(0) = (\tilde{D}/2C)v'(0)^2$, we see $u''(0) \not= 0$ if and only if 
$\tilde{D} \not= 0$. 

As seen above, the geometric solution $f$ is determined by 
two functions $p(u, v)  = p(v)$ and $y(u, v) + p_v(v)u$ of one variable $v$. 
From a pair of two functions $(\varphi(v), \psi(v))$ we get a geometric solution 
$$
f(u, v) = (u, \psi(v) - \varphi_v(v)u, \int (\varphi - v\varphi_v)du + v(\psi_v - u\varphi_{vv})dv, \varphi(v), v)
$$
to Hess $ = 0$. 
We define a mapping from a neighborhood of $(u_0, v_0)$ 
to $\R^4$ by 
$$
(u, v) \mapsto 
(\tilde{B}, \tilde{C}, C, \tilde{D}) = 
(\psi_v(v) - \varphi_{vv}(v)u, 
\psi_{vv}(v) - \varphi_{vvv}(v)u, \varphi_{vv}(v), 
\psi_{vvv}(v) - \varphi_{vvvv}(v)u). 
$$ 

For a generic $(\varphi(v), \psi(v))$, we have that 
$\varphi_{vv}$ vanishes 
at a finite number of points, where $\psi_v(v) \not= 0$. 
Therefore $(\tilde{B}, C) \not= (0, 0)$. 
Moreover two vectors 
$$
(\psi_v(v), \psi_{vv}(v), \psi_{vvv}(v)) \ 
{\rm{and}} 
\ 
(\varphi_{vv}(v), \varphi_{vvv}(v), \varphi_{vvvv}(v))
$$
are linearly independent 
for any $v$. 
Then, if $\tilde{B} = \psi_v(v) - \varphi_{vv}(v)u = 0$ for some 
$(u, v)$, then $(\tilde{C}, \tilde{D}) = 
(\psi_{vv}(v) - \varphi_{vvv}(v)u, \psi_{vvv}(v) - \varphi_{vvvv}(v)u) \not= (0, 0)$. 

Thus only the following cases occur 
on the $uv$-plane for a generic geometric solution: 

\ 

\begin{center} 
(i) $\tilde{B} \not= 0$, \quad (ii) $\tilde{B} = 0, \tilde{C} \not= 0, C \not= 0$, 
\quad (iii) $\tilde{B} = 0, \tilde{C} = 0, C \not= 0, \tilde{D} \not= 0$. 
\end{center}

\ 

In the case (i), $\pi\circ f$ is an immersion, in the 
case (ii), $\pi_1\circ f$ is a cuspidal edge, 
and, in the case (iii), $\pi_1\circ f$ is a swallowtail. 
Thus we have Theorem \ref {dev-front}.

\ 

To show Proposition \ref{OU}, we recall Lemma 11 of 
\cite{IMo}: 
\bel
\label{full}
Let $f : (\R^2,0) \to M$ be an open umbrella and $H : (M,f(0)) \to 
(\R,0)$ a smooth function-germ.  Suppose 
$dH(f(0)) \not= 0$ and $H\circ f = 0$. Then 
the kernel of $dH(f(0)) : T_{f(0)}M \to \R$ 
coincides with the contact hyperplane in $T_{f(0)}M$. 
\enl

An integral mapping $f : (\R^2, 0) \to M$ is called {\it full} if it has the property 
in Lemma \ref{full}. 
An integral mapping $f$ is called {\it formally full} 
if, for any formal complex valued function
$\hat{H}$ at $f(0)$ in $M$,  
the conditions $\hat{H}\circ \hat{f}_{\C} = 0$ and $d\hat{H}(0) \not= 0$ imply  
that the kernel of $d\hat{H}(f(0)) : T_{f(0)}M\otimes\C \to \C$ 
coincides with the complexification of the contact hyperplane in $T_{f(0)}M$. 
Here $\hat{f}_{\C}$ is the complexification of the formal Taylor 
series of $f$ at $0$ in $\R^2$. 
If $f$ is full, then it is formally full. 

For our treatment including the elliptic cases $(c > 0)$, 
we need a slightly strict result than Lemma \ref{full}: 
\bel
\label{formally full}
An open umbrella $f : (\R^2,0) \to M$ is formally full of rank one. 
\enl

\Proof 
By the contact invariance of the assumption and conclusion, we 
may suppose 
$$
f = (u, v^2, uv^3, v^3, \dfrac{3}{2}uv) : (\R^2, 0) \to (\R^5, 0)
$$ 
with respect to the contact form $\theta = dz - pdx - qdy$. 
Set 
$$
\hat{H}(x, y, z, p, q) = Ax + By + Cz + Dp + Eq + F(x, y, z, p, q), 
$$
where $A, B, C, D, E$ are complex numbers and 
$F$ is a formal power series starting from second order terms. 
The assumption means that 
$A, B, C, D, E$ are not all zero, and that 
$$
Au + Bv^2 + Cuv^3 + Dv^3 + E \dfrac{3}{2}uv 
+ F(u, v^2, uv^3, v^3, \dfrac{3}{2}uv) = 0, 
$$
as a formal power series of $u$ and $v$. 
Then we see $A = 0, B = 0, D = 0, E = 0$. 
Thus the kernel of $d\hat{H}(0)$ coincides with $\{\theta = dz = 0 \}$ 
at $0$ in $\R^5$. 
\QED

Then Proposition \ref{OU} follows from the following: 

\bel
\label{never formally full}
Let $f : (\R^2, 0) \lon M$ be an integral 
map-germ. Suppose $f$ is a generalized geometric solution 
to ${\mbox{\rm{Hess}}} = c, c \not= 0$, (resp. 
$K = c, c \not= 0$), 
of rank $\geq 1$. 
Then $f$ is never formally full.  
\enl

\Proof 
First we consider the case ${\rm{Hess}} = c$, $c \not= 0$. 
By the assumption we have 
$f^*\theta = 0, f^*\omega = 0$,  where 
$\theta = dz - pdx - qdy, \omega = c dx \wedge dy - dp \wedge dq$.
Now consider the complex valued $2$-forms
$$
\begin{array}{rcl}
(\sqrt{-c}dx + dq)\wedge (\sqrt{-c}dy - dp) 
& = & - cdx\wedge dy + dp \wedge dq + \sqrt{-c}(dp\wedge dx + dq\wedge dy) \\
& = & - \omega + \sqrt{-c}d\theta, 
\end{array}
$$
and 
$$
\begin{array}{rcl}
(\sqrt{-c}dx - dq)\wedge (\sqrt{-c}dy + dp) 
& = & - c dx\wedge dy + dp \wedge dq - \sqrt{-c}(dp\wedge dx + dq\wedge dy) \\
& = & - \omega - \sqrt{-c}d\theta. 
\end{array}
$$
Then we see $f^*(\sqrt{-c}dx + dq)\wedge f^*(\sqrt{-c}dy - dp) = 0$ 
and $f^*(\sqrt{-c}dx - dq)\wedge f^*(\sqrt{-c}dy + dp) = 0$, 
near $0$ on $\R^2$. 
On the other hand, 
since $f$ is of rank $\geq 1$, we see 
$$
(f^*dx)(0), (f^*dy)(0), (f^*dp)(0), (f^*dq)(0)
$$ are not all zero. 
Therefore 
$$
f^*(\sqrt{-c}dx + dq), \ f^*(\sqrt{-c}dy - dp), \ 
f^*(\sqrt{-c}dx - dq), \ f^*(\sqrt{-c}dy + dp)
$$ 
are not all zero at $0$ in $\R^2$. 
Suppose, for instance, $f^*(\sqrt{-c}dx + dq)(0) \not= 0$. 
Then, regarding $\sqrt{-c}x + q$ as a formal power series with complex coefficients, 
we have $(\sqrt{-c}x + q)\circ \hat{f}_{\C} = au + bv + \cdots$, 
with $(a, b) \not= (0, 0)$. 
Since $f^*(\sqrt{-c}dx + dq)\wedge f^*(\sqrt{-c}dy - dp) = 0$, 
we see  
there exists a formal power series $\hat{H}_1$ of 
$\sqrt{-c}x + q, \sqrt{-c}y - p$  
satisfying 
$\hat{H}_1((\sqrt{-c}x + q)\circ\hat{f}_{\C}, (\sqrt{-c}y - p)\circ\hat{f}_{\C}) = 0$ 
and 
$d\hat{H}(0) \not= 0$. 
Then we may set $\hat{H}(x, y, z, p, q) = 
\hat{H_1}(x + \sqrt{-c}q, y - \sqrt{-c}p)$. Therefore $f$ is not formally full. 

In the case $K = c$, $c \not= 0$, as seen in the previous section, we have 
$f^*\theta = 0, f^*\omega = 0$, where 
$$
\theta = dz - pdx - qdy, \quad 
\omega = c(1+p^2+q^2)^2dx\wedge dy - dp\wedge dq. 
$$
Since also in this case we have similar decompositions into two ways 
$$
(\sqrt{-c}(1+p^2+q^2)dx + dq)\wedge (\sqrt{-c}(1+p^2+q^2)dy - dp) 
= - \omega + \sqrt{-c}d\theta, 
$$
and 
$$
(\sqrt{-c}(1+p^2+q^2)dx - dq)\wedge (\sqrt{-c}(1+p^2+q^2)dy + dp) 
 =  - \omega - \sqrt{-c}d\theta, 
$$
we have a formal complex valued function $\hat{H}(x, y, z, p, q)$ with 
$\hat{H}\circ \hat{f}_{\C} = 0$, 
$d\hat{H}(0) \not= 0$, 
$\{ d\hat{H}(0) = 0 \} \not= D_{f(0)}$,  
the complexified contact hyperplane in $T_{f(0)}M\otimes\C$. 
Therefore $f$ is not formally full. 
\QED

\ber
{\rm 
In general, for a hyperbolic 
Monge-Amp\`{e}re system $\{ \theta, \omega \}$ on $\R^5$, 
there exist independent $1$-forms $\omega_1, \omega_2, \omega_3, \omega_4$ 
such that $\omega \equiv \omega_1\wedge \omega_2 \equiv 
\omega_3\wedge\omega_4$ mod.$d\theta$. 
Even in an elliptic case, we have similar decompositions into formal $1$-forms 
in the complex category. See for instance \cite{IL}\cite{BGH}. 
Thus the same proof as in Lemma \ref{never formally full} works to conclude the 
assertion: }
If an integral map-germ $f : (\R^2, 0) \to (\R^5, 0)$ 
is formally full and of rank $1$, then 
$f$ can not be a generalized geometric solution to any hyperbolic nor any elliptic  
Monge-Amp\`{e}re system on $\R^5$. 
\enr


\

\

\noindent Go-o ISHIKAWA \\
Department of Mathematics, Hokkaido University, 
Sapporo 060-0810, JAPAN. 
{\vspace{-0.3truecm}} 
\begin{verbatim}
E-mail : ishikawa@math.sci.hokudai.ac.jp
\end{verbatim}

\

\noindent Yoshinori MACHIDA \\ 
Numazu College of Technology, 
3600 Ooka, Numazu-shi, Shizuoka, 410-8501, JAPAN.  
{\vspace{-0.3truecm}} 
\begin{verbatim}
E-mail : machida@numazu-ct.ac.jp
\end{verbatim}

\end{document}